\documentclass[12pt]{article}
\usepackage{mathtools} 
\usepackage{amsmath, amssymb} 
\usepackage{amsthm,hyperref}
\usepackage{enumerate}  % lets you do \begin{enumerate}[(a)] and such
\usepackage{url} % formats url's nicely using \url{}

\usepackage{graphicx}

\DeclarePairedDelimiter\abs{\lvert}{\rvert}%
\DeclarePairedDelimiter\norm{\lVert}{\rVert}%

\makeatletter
\let\oldabs\abs
\def\abs{\@ifstar{\oldabs}{\oldabs*}}
\let\oldnorm\norm
\def\norm{\@ifstar{\oldnorm}{\oldnorm*}}
\makeatother

\newcommand\comp[1]{{\mkern2mu\overline{\mkern-2mu#1}}}

\newtheoremstyle{plainsl}%
	{\topsep}
	{\topsep}
	{\slshape} % only non-default setting
	{}
	{\normalfont\bfseries}
	{.}
	{ }
	{}

\swapnumbers
\numberwithin{equation}{section}

{\theoremstyle{plainsl}
\newtheorem{theorem}{Theorem}[section]
\newtheorem{lemma}[theorem]{Lemma}
\newtheorem{corollary}[theorem]{Corollary}}
{\theoremstyle{remark}
}

\renewcommand\proof{\noindent\textsl{Proof. }}
\newcommand\sqr[2]{{\vbox{\hrule height.#2pt
    \hbox{\vrule width.#2pt height#1pt \kern#1pt
        \vrule width.#2pt}\hrule height.#2pt}}}
\renewcommand\qed{%
	\ifmmode\eqno\sqr53
	\else\nolinebreak\ \hfill\sqr53\medbreak\fi}

\title{Equiangular Lines and Covers of the Complete Graph}
\author{G. Coutinho\footnote{Dep. of Combinatorics and Optimization, University of Waterloo. \newline \texttt{ \{gcoutinho, cgodsil, h3zhan\}@uwaterloo.ca}}, C. Godsil\footnotemark[1], M. Shirazi\footnote{\texttt{hamedshirazi@alumni.uwaterloo.ca}}, H. Zhan\footnotemark[1]}

\begin{document}
\maketitle

\begin{abstract}
The relation between equiangular sets of lines in the real space and distance-regular double covers of the complete graph is well known and studied since the work of Seidel and others in the 70's. The main topic of this paper is to continue the study on how complex equiangular lines relate to distance-regular covers of the complete graph with larger index. Given a set of equiangular lines meeting the relative (or Welch) bound, we show that if the entries of the corresponding Gram matrix are prime roots of unity, then these lines can be used to construct an antipodal distance-regular graph of diameter three. We also study in detail how the absolute (or Gerzon) bound for a set of equiangular lines can be used to derive bounds of the parameters of  abelian distance-regular covers of the complete graph. 
\end{abstract}

\section{Introduction}

We explore the rich relation between equiangular lines in a real or complex vector space and covers of the complete graph. In our journey, we will link these concepts to other structures and translate properties across the different topics.

Equiangular lines have been studied for a long a time, but there has been a recent surge in interest since the connection with quantum information theory was established (see, for example, Appleby \cite{appleby2005symmetric} or Scott and Grassl \cite{scott2010symmetric}). A set of equiangular lines meeting the so-called \textit{absolute (or Gerzon) bound} is also known as a \textit{symmetric, informationally complete, positive operator valued measured (SIC-POVM)}, and the problem of constructing SIC-POVMs is a major problem both in quantum information theory and in combinatorics. Also another important connection is the well known correspondence between equiangular lines and Seidel matrices, and that a set of lines meets the so-called \textit{relative (or Welch) bound} if and only if the corresponding Seidel matrix has only two distinct eigenvalues. These structures are also studied in frame theory, and sets of lines meeting the relative bound are called \textit{equiangular tight frames}.

In this paper, we study the relation between Seidel matrices and simple graphs that are covers of the complete graph. If a Seidel matrix has only two eigenvalues and if its entries are prime roots of unity, then we will show that such matrix implies the existence of distance-regular covers of the complete graph whose automorphism group satisfies certain properties, namely, that the covers are cyclic (in the sense defined by Godsil and Hensel \cite{Godsil1992}). This is a natural generalization of the well known correspondence between real tight frames and regular two-graphs.

Furthermore, we will use this relation to derive bounds on the defining parameters of such graphs using the absolute bound for equiangular lines. It turns out that the existence of certain graphs could give $d^2$ equiangular lines in $\mathbb{C}^d$. We give the parameter sets of these graphs, and analyse the real case similarly.

\section{Background on lines}

A set of lines spanned by vectors $x_1,...,x_n$ in $\mathbb{C}^d$ (or $\mathbb{R}^d$) is a set of \textit{complex (or real) equiangular lines} if there is $\alpha \in \mathbb{R}$ such that, for all $i$ and $j$,
\[|\langle x_i , x_j \rangle| = \alpha.\]
We will call $\alpha$ the \textit{angle} between two lines.

Upon associating each line determined by $x_i$ with the corresponding orthogonal projection given by $x_i x_i^*$, a set of equiangular lines is precisely the same thing as a $1$-regular quantum design of degree $1$. Such design will be called a \textit{tight frame} if
\[x_1x_1^{\ *} + ... + x_d x_d^{\ *} = \frac{n}{d} I.\]

\begin{theorem}[Relative bound, Van Lint and Seidel \cite{vanLintandSeidel}]
If there is a set of $n$ equiangular lines in dimension $d$, and if the angle of the set is $\alpha$, then
\[\alpha^2 \geq \frac{n-d}{(n-1) d}.\]
Equality holds if and only if the set of lines corresponds to a tight frame.
\end{theorem}

This result is sometimes called the Welch bound (see \cite{welch1974lower}).

Note that in the theorem above, it is irrelevant whether the set of lines is in a real or a complex space. The bound below however distinguishes between the two cases.

\begin{theorem}[Absolute bound, Gerzon (private communication to Seidel and Lemmens, see \cite{Lemmens1973})] \label{thm:absolutebound}
If there is a set of $n$ equiangular lines in $\mathbb{C}^d$, then
\[n \leq d^2.\]
If there is a set of $n$ equiangular lines in $\mathbb{R}^d$, then
\[n \leq \binom{d+1}{2}.\]
In either case, if equality holds, then the set of lines corresponds to a tight frame, and therefore the relative bound holds with equality.
\end{theorem}

Equiangular sets of lines are equivalent to other combinatorial structures, which we describe below. A \textit{Seidel matrix} is a Hermitian matrix whose diagonal entries are zeros and off-diagonal entries have absolute value one.

\begin{theorem}[Lemmens and Seidel \cite{Lemmens1973}]
Let $S$ be an $n\times n$ Seidel matrix. Let $m_\tau$ be the multiplicity of its least eigenvalue $\tau$, and $m_\theta$ the multiplicity of its largest eigenvalue $\theta$. Then 
\begin{itemize}
\item the matrix $I - (1/\tau) S$ is the Gram matrix of a set of $n$ equiangular lines in dimension $n-m_\tau$, and
\item the matrix $I - (1/\theta) S $ is the Gram matrix of a set of $n$ equiangular lines in dimension $n-m_\theta$.
\end{itemize}
Moreover, a set of $n$ equiangular lines over $\mathbb{C}^d$ with angle $\alpha$ meets the equality in the relative bound if and only if the corresponding Seidel matrix (up to a sign) has precisely two distinct eigenvalues, and they are
\[-\frac{1}{\alpha} \ , \quad \frac{n-d}{\alpha d}\]
with multiplicities $n-d$ and $d$, respectively.
\label{Seidel}
\end{theorem}

\section{Background on covers}

For more details about the content of this section, we refer the reader to Godsil and Hensel \cite{Godsil1992}.

An \textit{$r$-fold cover of $K_n$} is a simple graph $X$ on $r n$ vertices satisfying the following two properties:
\begin{enumerate}[(1)]
\item There is a partition of the vertex set of $X$ into $n$ sets of $r$ vertices each, to be called \textit{fibres}, such that no two vertices in the same fibre are connected.
\item There is a perfect matching between any two fibres.
\end{enumerate}
\newcommand{\drackn}{\textsc{drackn}}
In some cases, an $r$-fold cover of $K_n$ will enjoy the property of being distance-regular. When this happens, the distance-regular graph will be antipodal and all of its intersection parameters will be determined by $n$, $r$ and a third parameter, typically denoted by $c$, which counts the number of common neighbours of two vertices at distance two. We will shortly refer to such a graph as an $(n,r,c)$-\drackn \ (standing for \textit{distance-regular antipodal cover of $K_n$}). Conversely, one can show that all antipodal distance-regular graphs of a diameter $3$ are covers of a complete graph. It is straightforward to show that $n-1$ and $-1$ are eigenvalues of these graphs with respective multiplicities $1$ and $n-1$. There are another two distinct eigenvalues, of opposing signs, that are typically called $\theta$ and $\tau$ with the convention that $\theta > 0$ and $\tau < 0$. If $\delta = n - rc - 2$, it follows that
\begin{align}
\theta = \frac{\delta + \sqrt{\delta^2 + 4(n-1)}}{2} \quad \text{and} \quad \tau = \frac{\delta - \sqrt{\delta^2 + 4(n-1)}}{2} \label{eq:eigenvalues}
\end{align}

The problem of determining which parameter sets correspond to an actual distance-regular graph has been attacked for decades, and despite the many efforts, its full solution most likely will not be seen in the near future. This is the case even for distance-regular graphs of small diameter and constrictive structural properties, such as \drackn s. However, we can compile a list of fairly restrictive non-trivial conditions that $n$, $r$ and $c$ must satisfy in order to correspond to a graph.

\begin{theorem}
Let $X$ be an $(n,r,c)$-\drackn \ with distinct eigenvalues $n-1$,  $\theta$, $-1$ and $\tau$, and with $n\ge 2, r\ge 2, c\ge 1$. Recall that $\delta = n - rc -2$. Then the parameters of $X$ satisfy the following conditions.
\begin{enumerate}[(a)]
\item $1\le c(r-1)\le n-2 \leq c(2r - 1) - 2$.
\item The multiplicities of $\theta$ and $\tau$ satisfy
\[m_{\theta}=\frac{n(r-1)\tau}{r-\theta} \quad \text{and} \quad m_{\tau}=\frac{n(r-1)\theta}{\theta-\tau},\]
and these ratios must be integers.
\item If $\delta \neq 0$, then $\theta$ and $\tau$ are integers.
\item If $\delta = 0$, then $\theta=-\tau=\sqrt{n-1}$. 
\item If $n$ is even, then $c$ is even.
\item If $c=1$, then $(n-r)$ divides $n-1$,  $(n-r)(n-r+1)$ divides $rn(n-1)$, and $(n-r)^2\le n-1$.
\item If $r>2$, then $\theta^3\ge n-1$.
\item Suppose $\theta\ne 1$, $\tau\ne -1$, and $\theta^3\ne n-1$. For $r>2$, we have
\[rn\le \frac{1}{2}m_{\theta}\left(m_{\theta}+1\right),\quad rn\le \frac{1}{2}m_{\tau}\left(m_{\tau}+1\right).\]
For $r=2$, we have
\[n\le \frac{1}{2}m_{\theta}\left(m_{\theta}+1\right),\quad n\le \frac{1}{2}m_{\tau}\left(m_{\tau}+1\right).\]
\item Let $r>2$ and $\beta\in\{\theta,\tau\}$ be an integer. If $n>m_{\beta}-r+3$, then $\beta+1$ divides $c$.
\end{enumerate}
\label{feascond}
\end{theorem}

An \textit{arc function} of index $r$ over $K_n$ is a function $f$ from the arcs of $K_n$ to the symmetric group $\textrm{Sym}(r)$ satisfying $f(u,v)^{-1} = f(v,u)$.
Arc functions are equivalent to covers in a very natural way. In an $r$-fold cover of $K_n$, a matching from the fibre corresponding to a vertex $u$ of $K_n$ to the fibre corresponding to a vertex $v$ of $K_n$ can be seen as a permutation on $r$ elements that is precisely equal to $f(u,v)$. Without loss of generality, we can always suppose that an arc function $f$ will be equal to the identity permutation when evaluated over the edges of a spanning tree. When this happens, $f$ is called a \textit{normalized arc function}.

Let $\langle f \rangle$ be the permutation group generated by the images of $f$ over all arcs of $K_n$. An $r$-fold cover of $K_n$ determined by a normalized arc function $f$ is called \textit{regular} if $\langle f \rangle$ is regular, and moreover \textit{abelian} if $\langle f  \rangle$ is an abelian group. If $\langle f  \rangle$ is a cyclic group, we say that the cover is \textit{cyclic}. Note that a cover is regular if and only if $|\langle f  \rangle| = r$. Finally, the automorphism group of the cover that fixes each fibre as a set is regular if and only if $\langle f \rangle$ is regular, and in this case these are isomorphic groups.

Consider the square matrix whose rows and columns are indexed by the vertices of $K_n$, and whose entry $(u,v)$ is equal to the permutation $f(u,v)$ and diagonal entries are equal to $0$. We denote this matrix by $A(K_n)^f$.

Let $\varphi$ be an $s$-dimensional representation of $\langle f \rangle$. Let $A(K_n)^{\varphi(f)}$ stand for the matrix obtained from $A(K_n)^f$ by replacing each of its entries by an $s \times s$ permutation matrix corresponding to their image under $\phi$, and the diagonal entries by $s \times s$ blocks of $0$.

If $X$ is a regular cover defined by $f$ and $\phi$ is the regular representation of $\langle f \rangle$, then
\[A(K_n)^{\phi(f)} = A(X).\]
We summarize in the next theorem some important facts relating covers, representations and linear algebra. This and more can be found in Godsil and Hensel \cite[Sections 8 and 9]{Godsil1992}. 

The first part is an immediate consequence of the well known expression for the eigenvectors of the regular representation of abelian groups in terms of linear characters. The second is a restatement of \cite[Corollary 7.5]{Godsil1992}.

\begin{theorem}\label{thm:characters}
Let $X$ be a connected abelian $r$-fold cover of $K_n$ determined by a normalized arc function $f$. Let $\phi_1,...,\phi_r$ be the linear characters of $\langle f \rangle$. Then $A(X)$ is similar to the matrix
\[\begin{pmatrix}
A(K_n)^{\phi_1(f)} & & & & \\  & A(K_n)^{\phi_2(f)} &  &  &  \\  &  &  \ddots &  &  \\  &  &  & & A(K_n)^{\phi_r(f)}
\end{pmatrix}.\]
Moreover, $X$ is an $(n,r,c)$-\drackn \ if and only if the minimal polynomial of each matrix in ${\{ A(K_n)^{\phi_i(f)} : i = 2,...,r \}}$ is
\[x^2 - (n-rc-2) x - (n-1).\]
\end{theorem}

\section{Lines from covers} \label{Lines from Covers}

A consequence of Theorem \ref{thm:characters} is that the existence of abelian \drackn s implies the existence of sets of complex equiangular lines meeting the relative bound.

\begin{theorem}\label{cor:covergivinglines}
Let $X$ be an abelian $(n,r,c)$-\drackn \ defined by a symmetric arc function $f$. Suppose the eigenvalues of $X$ are $n-1$, $\theta$, $-1$ and $\tau$, with respective multiplicities $1$, $m_\theta$, $n-1$ and $m_\tau$. Let $\phi$ be a non-trivial character of $\langle f \rangle$. Then $A(K_n)^{\phi(f)}$ is a Seidel matrix with precisely two distinct eigenvalues, $\theta$ and $\tau$, and therefore:
\begin{enumerate}[(a)]
\item There are $n$ complex equiangular lines in dimension $\left( n - \dfrac{m_\theta}{r-1} \right)$ meeting the relative bound.
\item There are $n$ complex equiangular lines in dimension $\left( n - \dfrac{m_\tau}{r-1} \right)$ meeting the relative bound.
\end{enumerate}
\end{theorem}
\proof
Let $\phi_1,...,\phi_r$ be the linear characters of $\langle f \rangle$. If $\phi_1$ is the trivial character, then $A(K_n)^{\phi_1(f)} = A(K_n)$, therefore its spectrum is $(n-1)^{(1)}$ and $(-1)^{(n-1)}$. By Theorem \ref{thm:characters}, this implies that the spectrum of each matrix $A(K_n)^{\phi_k(f)}$, for $k=2,...,r$, contains only multiple copies of $\theta$ or $\tau$. Because these matrices have trace $0$, the multiplicities of $\theta$ and $\tau$ in each one of them do not depend on $k$, and will be equal to $m_\theta / (r-1)$ and $m_\tau / (r-1)$ respectively. The result now follows from Theorem \ref{Seidel}.
\qed

We present here some examples of infinite families of abelian covers.

\subsubsection*{Symplectic covers}

This construction generalizes the so-called Thas - Somma construction. Let $p$ be a prime, and let $V$ and $U$ be vectors spaces over $GF(p)$ of respective dimensions $m$ and $s$. Let $B$ be a $GF(p)$-linear alternating form from $V \times V$ to $U$, with the extra property that for each $a \in V$, the linear mapping $B_a : V \to U$ defined as $B_a(v) = B(a,v)$ is surjective. We define a graph $X(B)$ on the vertex set $V \times U$ where adjacency between distinct vertices is defined by
\[(v,a) \sim (w,b) \iff B(v,w) = a - b.\]
It follows that $X(B)$ is a $(p^m,p^s,p^{m-s})$-\drackn . If $f$ is the normalized arc function defining $X(B)$, it is easy to see that $\langle f \rangle \cong \mathbb{Z}_p^{\ s}$, hence $X(B)$ is abelian. Godsil \cite{MR1414466} describes a way of constructing such symplectic forms whenever there exists an $s$-dimensional space of invertible $m \times m$ skew-symmetric matrices over $GF(p)$. In particular, this always exists when $s = 1$ and $m$ is even.

\subsubsection*{De Caen and Fon-der-Flaass construction}

Let $V$ be a $d$-dimensional vector space over $GF(2^t)$. Consider a skew product $*$ of $V$, that is, a bilinear mapping $*: V^2 \to V$ such that $x \mapsto x* x$ is a bijection and $x * y = y * x$ if and only if $x$ and $y$ are linearly dependent. Let $S = ((s_{ij}))_{i,j \in F}$ be a symmetric latin square filled with the elements of $GF(2^t)$. Construct a graph on the vertex set ${V \times F \times V}$ and such that $(a,i,\alpha)$ and $(b,j,\beta)$ are adjacent if and only if
\[\alpha + \beta = a * b + b* a + s_{ij}(a* a + b*b).\]

This graph is a $(2^{t(d+1)},\ 2^{td},\ 2^t)$-\drackn , and by construction it also follows that it is an abelian cover whose automorphism group fixing each fibre is isomorphic to $\mathbb{Z}_2^{\ td}$. De Caen and Fon-der-Flaas \cite{deCaenFonderFlaass} also remark that skew products exist if and only if $d$ is odd.

\subsubsection*{Generalized Hadamard matrices}

Let $X$ be an abelian $(n,r,c)$-\drackn \ determined by a normalized arc function $f$ such that $\delta = -2$. If $\phi$ is a character of $\langle f \rangle$, let $S = A(K_n)^{\phi(f)}$ (recall the notation introduced after Theorem \ref{feascond}). It follows from Theorem \ref{thm:characters} that
\[S^2 = (n-1)I + \delta S,\]
and so because $\delta = -2$, we have that $(S+I)^2 = n I$. The matrix $S+I$ is therefore an Hermitian Butson-type Hadamard matrix with constant diagonal. 

Let $G = \langle f \rangle$ and $e$ be the identity of $G$. Consider the group ring $\mathbb{Z}[G]$. Given a subset $S$ of $G$, we use the notation
\[\underline{S} = \sum_{g \in S} g,\]
where the sum is the ring sum. Consider the matrix $H = A(K_n)^f + e I$ over $\mathbb{Z}[G]$. Note that
\[H^2 = n I + c~\underline{G}~(J-I),\]
therefore $H$ is a \textit{generalized Hadamard matrix} $GH(r,c)$ (note that $n = rc$) over the group $G$ (of order $r$). See Colbourn and Dinitz \cite[Chapter V.5]{HandbookDesigns} for more details.

The paragraph above shows that any abelian $(rc,r,c)$-\drackn \ implies the existence of a self-adjoint $GH(r,c)$ with constant diagonal, and \cite[Corollary 7.5]{Godsil1992} states precisely the opposite.

In a recent paper, Klin and Pech \cite{klin2011new} derived many new constructions of abelian \drackn s based on these generalized Hadamard matrices. In particular, they showed in \cite[Theorem 5.6]{klin2011new} that any $n \times n$ generalized Hadamard matrix over a group $G$ implies the existence of a $n^2 \times n^2$ self-adjoint generalized Hadamard matrix with constant diagonal over $G$, and by the remarks above, those are equivalent to abelian \drackn s with $\delta = -2$. The table below, extracted from \cite[page 227]{klin2011new}, contains a list of the parameter sets of \drackn s that are obtained from known generalized Hadamard matrices.

\renewcommand{\arraystretch}{1.6}
\[\begin{array}{c|c}
 (n,r,c) & \text{conditions} \\\hline
(p^{m 2^t},p^{n},p^{m 2^t - n}) & m \geq n \geq 1, \text{ $p$ prime, } t > 0 \\  \hline
(2^{2^t} p^{m 2^t},p^{n},p^{m 2^t - n}) &  m \geq n \geq 1, \text{ $p$ prime, } t > 0 \\ \hline
(4^{2^t} p^{m 2^t},p^{n},4^{2^t} p^{m 2^t - n}) &  m \geq n \geq 1, \text{ $p$ prime, } t > 0 \\  \hline
(8^{2^t} q^{2^t},q,8^{2^t}p^{2^t - 1}) & 19 < q < 200, \text{ $q$ prime power, } t > 0   \\  \hline
(8^{2^t} p^{2^t},p,8^{2^t}p^{2^t - 1}) &  19 < p, \text{ $p$ prime, } t > 0 \\  \hline
(k^{2^t} q^{2^t},q,k^{2^t} q^{2^t - 1}) &  q > ((k-2)2^{k-2} ), \ \exists \text{ Hadamard matrix of order $k$, } t > 0 \\ \hline
(45,3,12) & \text{ due to Klin and Pech, coming from the Foster graph. } \\  \hline
(144,4,36) & \text{ discovered by J. Seberry.} 
\end{array}\]

\section{Covers from lines} \label{section covers from lines}

In Section \ref{Lines from Covers}, we showed how to construct equiangular lines using an abelian \drackn. On the opposite direction, it is well known that real equiangular lines can be used to construct regular two-graphs, as we will briefly explain. Then, we generalize this result, showing how a set of complex equiangular lines can be used to construct an abelian \drackn\ of index larger than two.

Any given $(n,2,c)$-\drackn\ is necessarily abelian. These graphs are equivalent to the so-called and well studied regular two-graphs. The Seidel matrix $S$ of any set of real equiangular lines has off diagonal entries equal to $+1$ or $-1$, and upon replacing
\[+1 \quad \text{by} \quad \begin{pmatrix}
1 & 0 \\ 0 & 1
\end{pmatrix} \qquad \text{and} \qquad -1 \quad \text{by} \quad \begin{pmatrix}
0 & 1 \\ 1 & 0
\end{pmatrix},\]
it is easy to see that the new $2n \times 2n$ matrix is the adjacency matrix of a $2$-fold cover of $K_n$, say $X$. If the relative bound is met on the original set of lines, then $S$ has only two distinct eigenvalues, and so $X$ will be an abelian \drackn. This describes the correspondence between real equiangular tight frames and regular two-graphs. The classification of regular two-graphs has received a considerable amount of attention in the past 40 years or so. See Godsil and Royle \cite[Chapter 11]{Godsil2001} or Brouwer and Haemers \cite[Chapter 10]{BrouwerHaemers} for more information.

We will now describe a way of constructing cyclic \drackn s of index larger than two from sets of equiangular lines satisfying certain properties.

\begin{theorem} \label{thm:coversfromlines}
Suppose $x_1,...,x_n$ is a set of complex equiangular lines in $\mathbb{C}^d$ with angle $\alpha$ and Gram matrix $G$. Suppose they satisfy the following two properties:
\begin{enumerate}[(i)]
\item This set of lines meets the relative bound, and hence
\[\alpha^2 = \frac{n-d}{(n-1)d}.\]
\item All off-diagonal entries of the matrix
\[S = \frac{1}{\alpha} (G-I)\]
are $r$-th roots of unity, where $r$ is a prime.
\end{enumerate}
Then there exists a cyclic $(n,r,c)$-\drackn , where
\[c =\frac{1}{r} \left( (n-2) + \frac{2d-n}{\alpha d} \right).\]
\end{theorem}

\proof
Let $C_r$ be the multiplicative group of the $r$-th roots of unity. Let $\varphi$ be a representation of $C_r$ of degree $k$, and let $S^\varphi$ be the matrix obtained from $S$ by replacing each diagonal entry by a $k \times k$ block of $0$s, and each off-diagonal entry by its image under $\varphi$. Naturally, if $\phi$ is the regular representation of $C_r$, $S^\phi$ is the adjacency matrix of a graph $X$, and our goal is to show that $X$ is a \drackn .

Let $\phi_1,...,\phi_r$ be the linear characters of $C_r$, satisfying
\[\phi_k (\mathrm{e}^{2 \pi \mathrm{i} /r}) = \mathrm{e}^{(k-1)2 \pi \mathrm{i} /r}.\]
Clearly $A(K_n) = S^{\phi_1}$ and $S = S^{\phi_2}$. We claim that, for $j = 2,...,r$, the minimal polynomials of the matrices $S^{\phi_j}$ are all equal.

In fact, let $\Phi_r(x)$ be the $r$-th cyclotomic polynomial, and because $r$ is prime, we have ${\Phi(x) =  x^{r-1} + x^{r-2} + ... + x + 1}$. Since
\[\mathbb{Q}(\mathrm{e}^{2 \pi \mathrm{i} / r}) \cong \mathbb{Q}[x] / \langle \Phi_r(x) \rangle,\]
we can see the matrix $S$ as a matrix whose off-diagonal entries are powers of the indeterminate $x$ subject to the relation $\Phi(x) = 0$. Thus, if $m(y)$ is the minimal polynomial of $S$, then $m(S^{\phi_j})$ will vanish because $\phi_j(\mathrm{e}^{2 \pi \mathrm{i}  / r})$ is also a root of $\Phi_r(x)$ for all $j$. From Theorem \ref{Seidel} and condition (i), $m(y)$ has degree two and $S^{\phi_j}$ is not a multiple of the identity matrix for any $j$, hence the claim follows.

Moreover, the trace of $S^{\phi_j}$ is equal to $0$ for all $j$, therefore all these matrices with $j \geq 2$ are cospectral. Because the eigenvectors of the regular representation of an abelian group are its characters, it follows that $S^\phi$ is similar to the block diagonal matrix whose blocks are the matrices $S^{\phi_j}$, with $j = 1,...,r$. All together, and from the expression for the eigenvalues of $S$ given by Theorem \ref{Seidel}, we have shown that $X$ is an $r$-fold cover of $K_n$ with spectrum given by
\[n-1^{(1)}, \quad \left(\frac{n-d}{\alpha d}\right)^{\big(d (r-1) \big)}, \quad -1^{(n-1)} ,\quad  \left(\frac{-1}{\alpha}\right)^{\big( (r-1)(n-d) \big)}.\]
By Godsil and Hensel \cite[Lemma 7.1]{Godsil1992} and because $r$ is prime, it follows that $X$ is connected.

Finally, by Theorem \ref{thm:characters}, we have that $X$ is a \drackn \ with the given parameters.
\qed

Unfortunately, the only examples we know of sets of lines satisfying the conditions of Theorem \ref{thm:coversfromlines} are those constructed from known abelian \drackn s. Recently, Fickus et al.~(see \cite{fickus2012steiner} and \cite{fickuskirkman}) showed how to construct sets of equiangular lines meeting the relative bound based on previous sets and combinatorial designs. Using their construction, one is almost capable of obtaining a new set of lines satisfying the conditions of Theorem \ref{thm:coversfromlines}. The only problem is that certain entries of the matrix $S$ will be equal to $-1$, and thus not an $r$-th root of unity for any $r$ prime other than $r=2$. In this case, the corresponding abelian \drackn s coming from certain designs were already known (see Goethals and Seidel \cite{goethalsseidelSRGsandDesigns}).

\section{Feasibility conditions for DRACKNs}

In Theorem \ref{feascond}, we presented feasibility conditions for the parameter sets of \drackn s. In this section, we work out some extra feasibility conditions that the parameters of an abelian \drackn \ must satisfy. Our main tool will be to use the absolute bound for a set of equiangular lines to find bounds on the parameters of abelian \drackn s. We will then study the extreme cases. For instance, we find, to our surprise, that there are some feasible parameter sets of abelian \drackn s that would give $d^2$ equiangular lines in $\mathbb{C}^d$.

We begin by pointing out an immediate corollary of two results due to Godsil and Hensel.

\begin{theorem}[``quotienting", \cite{Godsil1992}, Lemma 6.2]\label{quotienting}
Let $X$ be an abelian $(n,r,c)$-\drackn \ determined by a normalized arc function $f$. Let $H$ be a subgroup of $\langle f \rangle$ of size $t$. Then the partition induced by the orbits of $H$ in each fibre is equitable, and therefore there is an abelian $(n,r/t,tc)$-\drackn \ obtained as a quotient by this partition.
\end{theorem}
\begin{theorem}[\cite{Godsil1992}, Theorem 9.2]
If $X$ is a cyclic $(n,r,c)$-\drackn \ with $r>2$, then $r$ divides $n$.
\end{theorem}
If $p$ is a prime that divides the order of a group $G$, then there is a cyclic subgroup of $H$ of $G$ of order $p$. Using this, we obtain the corollary below, which in particular implies that \drackn s with $\delta = 0$ and $r$ not a power of two cannot be abelian.
\begin{corollary}\label{cor:r divides n}
If $X$ is an abelian $(n,r,c)$-\drackn , then any odd prime that divides $r$ also divides $n$.
\end{corollary}

Now we show how to translate the absolute bound for a set of equiangular lines into some extra feasibility conditions for the parameter sets of abelian \drackn s. To simplify the notation, when $r$ is given by the context and $m$ is an integer, let
\[\overline{m} = \frac{m}{r-1}.\]
Throughout the following results, the parity of $r$ will play an important role. The reason is that every abelian group of even order has a real-valued linear character. The corresponding Seidel matrix of this character is a matrix with real entries, and hence corresponds to a set of real equiangular lines. Hence if $r$ is even, the parameters of the \drackn\ will be subject to the absolute bound for real lines, and therefore the bounds will be more restrictive.

\begin{lemma}
Let $X$ be an abelian $(n,r,c)$-\drackn \ with distinct eigenvalues $n-1 > \theta > -1 >\tau$. If $r$ is even, then
\[-\frac{1}{2}\sqrt{(n-1)\left(\sqrt{8n+1}-3\right)}\le \tau \le -\sqrt{\frac{1}{2}\left(\sqrt{8n+1}+3\right)}.\]
Moreover, if the lower bound is tight, then there is a set of real equiangular lines of size $\binom{\comp{m_{\tau}}+1}{2}$ in dimension $\comp{m_{\tau}}$, and if the upper bound is tight, then there is a set of real equiangular lines of size $\binom{\comp{m_{\theta}}+1}{2}$ in dimension $\comp{m_{\theta}}$.

If $r$ is odd, then
\[-\left(\sqrt{n}-1\right)\sqrt{\sqrt{n}+1} \le \tau\le -\sqrt{\sqrt{n}+1}.\]
Moreover, if the lower bound is tight, then there is a set of complex equiangular lines of size $\overline{m_\tau}^{\hspace{0.3mm} 2}$ in dimension $\overline{m_\tau}$, and if the upper bound is tight, then there is a set of complex equiangular lines of size $\overline{m_\theta}^{\hspace{0.3mm} 2}$ in dimension $\overline{m_\theta}$.
\label{lem:bounds}
\end{lemma}
\proof
We explicitly prove the second lower bound. The other three bounds are very similar. By Theorem \ref{cor:covergivinglines}, the existence of the \drackn\ implies the existence of $n$ complex equiangular lines in dimension $\overline{m_\tau}$. By the absolute bound \ref{thm:absolutebound},
\[n\leq \overline{m_{\tau}}^{\hspace{0.3mm} 2}.\]
From the expression for the multiplicities given by Theorem \ref{feascond},
\[\overline{m_{\tau}}=\frac{n}{1+\tau^2/(n-1)}\]
which is strictly decreasing in $\abs{\tau}$, and is greater than or equal to $\sqrt{n}$ if and only if
\[\abs{\tau} \leq \left(\sqrt{n}-1\right)\sqrt{\sqrt{n}+1}.\]
Since $\tau<0$, the lower bound follows. 
\qed
Below, we show that the extreme cases in the bounds of the lemma above can be conveniently parametrized.
\begin{theorem}
Let $X$ be an abelian $(n,r,c)$-\drackn\ with $r$ even. It gives a set of real equiangular lines meeting the absolute bound if and only if, for some positive integer $t$ or $t = \sqrt{5}$, one of the following cases holds.

\[\begin{array}{c|c|c}
\text{Parameter} & \text{case (I.a)} & \text{case (I.b)} \\\hline
n & \frac{1}{2}(t^2-2)(t^2-1) & \frac{1}{2}(t^2-2)(t^2-1) \\  \hline
rc & \frac{1}{2}(t+1)^3(t-2) & \frac{1}{2}(t-1)^3(t+2)\\ \hline
\delta & -\frac{1}{2}t(t^2-5) & \frac{1}{2}t(t^2-5) \\  \hline
\theta & t & \frac{1}{2}t(t^2-3)\\  \hline
\tau & -\frac{1}{2}t(t^2-3) & -t \\  \hline
\overline{m_\theta} &  \frac{1}{2}(t^2-2)(t^2-3)  & t^2-2 \\ \hline
\overline{m_\tau} & t^2-2 & \frac{1}{2}(t^2-2)(t^2-3) \\  
\end{array}\]
Let $X$ be an abelian $(n,r,c)$-\drackn\ with $r$ odd. It gives a set of complex equiangular lines meeting the absolute bound if and only if, for some positive integer $t$, one of the following cases holds.
\[\begin{array}{c|c|c}
\text{Parameter} & \text{case (II.a)} & \text{case (II.b)} \\ \hline
n & (t^2-1)^2 & (t^2-1)^2 \\ \hline
rc & (t+1)^2(t^2-t-1) & (t-1)^2(t^2+t-1)\\ \hline
\delta & -(t^2-3)t & (t^2-3)t \\ \hline
\theta & t & (t^2-2)t \\ \hline
\tau & - (t^2-2)t & -t \\ \hline
\overline{m_\theta} &  (t^2-2)(t^2-1)  & t^2-1\\ \hline
\overline{m_\tau} & t^2-1 & (t^2-2)(t^2-1) \\ 
\end{array}\]
\renewcommand{\arraystretch}{1
}
\label{thm:parameters}
\end{theorem}
\proof
It is trivial to check that if either of the cases hold, then the absolute bound is satisfied with equality. In the cases (a), we have a maximum sized set of lines in dimension $\overline{m_\tau}$, whereas in cases (b), the dimension of the lines is given by $\overline{m_\theta}$.

For the converse, note that each case corresponds to equality being achieved in either the lower or the upper bound in either of the cases of Lemma \ref{lem:bounds}. We show the case (I.a) of the converse. The other three cases are similar. Suppose the lower bound in the case where $r$ is even in Lemma \ref{lem:bounds} is tight. Thus
\begin{equation}
\tau=-\frac{1}{2}\sqrt{(n-1)\left(\sqrt{8n+1}-3\right)}. \label{tau}
\end{equation}
From (\ref{eq:eigenvalues}), it follows that $\delta = \tau - (n-1)/\tau$. If $\delta = 0$, then $\tau = - \sqrt{n-1}$. Plugging it into (\ref{tau}) gives $n=6$. The remaining parameters can be computed accordingly and fit the expressions in case (I.a) with $t = \sqrt{5}$. If $\delta\ne 0$, then by Theorem \ref{feascond}, both eigenvalues $\theta$ and $\tau$ are integers. Let $t = \theta \in \mathbb{Z}^+$. From (\ref{eq:eigenvalues}), we have $\theta \tau = - (n-1)$. Coupling with (\ref{tau}), it follows that
\[n=\frac{1}{2}(t^2-2)(t^2-1).\]
The remaining parameters can be calculated accordingly.
\qed

Now we proceed to unfold various cases presented in Theorem \ref{thm:parameters}. Our final goal is to present a complete list of parameter sets which satisfy all known feasibility conditions given by Theorem \ref{feascond} and Corollary \ref{cor:r divides n}. For the reasons presented in the beginning of Section \ref{section covers from lines}, we will skip the case $r = 2$.

\subsection{Case I.a}

\begin{theorem}
If $r\ge 4$, then the only feasible parameter set $(n,r,c)$ for an abelian \drackn\ corresponding to set of real equiangular lines in dimension $\comp{m_{\tau}}$ of maximum size is the $(28,4,8)$.
\end{theorem}

\proof
Suppose there is a $t$ such that $(n,r,c)$ and the other parameters are given as in case (I.a). Because $rc > 0$, $t > 2$. If $t = \sqrt{5}$, then in view of Theorem \ref{feascond} condition (f), the only possible parameter set is $(6,2,2)$. So suppose $t \geq 3$. By condition (h) in Theorem \ref{feascond}, 
\[t^3=\theta^3 \geq n-1=\frac{1}{2}(t^2-2)(t^2-1)-1,\]
which gives $t\geq 3$. Therefore $t=3$. The corresponding parameters are
\[n=28,\quad rc=32,\quad \theta=3,\quad \tau=-9.\]
Thus $r \in \{4,8,16,32\}$. If $r>4$, then $c \leq 4$, and so condition (a) of Theorem \ref{feascond} is not satisfied. Therefore $r=4$.
\qed

\subsection{Case II.a}

When $t = 2$, an abelian $(9,3,3)$-\drackn \ corresponding to case (II.a) above exists. It can be obtained via the Thas - Somma construction described in the end of Section \ref{Lines from Covers}. %TODO ATTENTION HERE 
The other cases are covered by the result below.

\begin{corollary}
The only parameter set $(n,r,c)$ corresponding to an abelian \drackn \ that gives a set of complex equiangular lines in dimension $\comp{m_{\tau}}$ of maximum size is the $(9,3,3)$.
\end{corollary}
\proof
By Theorem \ref{feascond}, if $r>2$, then 
\[\theta^3\ge n-1.\]
Plugging in $n=(t^2-1)^2$ and $\theta=t$ yields
\[t^2-t-2\le 0,\]
which restricts $t$ to 2.
\qed

\subsection{Case I.b}

\begin{theorem} \label{thm:caseIb}
Suppose that for some positive integer $t$, parameters $(n,r,c)$ satisfy
\[n = \frac{1}{2}(t^2 - 2)(t^2 - 1) \quad \text{and} \quad rc = \frac{1}{2}(t-1)^3 (t+2),\]
and $r$ is an even integer. These parameters satisfy all the feasibility conditions listed in Theorem \ref{feascond} and in Corollary \ref{cor:r divides n} for the existence of a corresponding abelian \drackn \ if and only if all of the following hold.
\begin{enumerate}[(1)]
\item $t\ge 3$ and is not divisible by four.
\item
$c\ge 2$.
\item
If $r \leq (1/2)(t^2+1)$, then $r$ divides $t-1$.
\item 
If $t$ is odd then $c$ is even.
\item
Any odd prime that divides $r$ must divide $t-1$ as well.
\end{enumerate}
\end{theorem}
\proof
First we show that conditions (1) - (5) are necessary.
\begin{itemize}
\item[(1)] If $t=2$, then $n=3$ and $r=2$. Hence $t\ge 3$. On the other hand, $r c$ is even if and only if $t$ is not divisible by four.
\item[(2)]
Assume that $c=1$. Then condition (f) in Theorem \ref{feascond} says that
\[(n-rc)^2=(n-r)^2\leq n-1.\]
However
\[4\left((n-r)^2-(n-1)\right)=((t^2-2)(t^2-1)-(t-1)^3(t+2))^2-2(t^2-2)(t^2-1)+4\]
is positive for all $t\ge 3$, contradicting our assumption. Thus $c\ge 2$.
\item[(3)]
With the given parameters, note that
\[n>m_{\theta}-r+3\]
and $t \geq 3$ if and only if
\[r\le \frac{1}{2}(t^2+1).\]
It also follows from the given parameters that
\[(\theta + 1) (t-1) = rc,\]
so $(\theta+1)$ divides $c$ if and only if $r$ divides $(t-1)$. Thus condition (i) of Theorem \ref{feascond} implies (3).
\item[(4)] 
Since 
\[n=\frac{1}{2}(t^2-2)(t^2-1)\]
we see that $n$ is even if and only if $t$ is odd. By condition (e) in Theorem \ref{feascond}, if $t$ is odd, then $n$ is even and so $c$ is even. 
\item[(5)]
This is the statement of Corollary \ref{cor:r divides n}.
\end{itemize}
Now we proceed to show that condition (1) - (4) are sufficient to guarantee all conditions in Theorem \ref{feascond}.
\begin{itemize}
\item[(a)] The first bound is satisfied as $r\ge 2$ and $c\ge 2$. The second is equivalent to $(\delta+c) \ge 0$, which is also true since
\[\delta=\frac{1}{2}t(t^2-5)\]
is positive whenever $t\ge 3$. The third bound is equivalent to
\[r\ge \frac{rc }{2rc -n}.\]
From the given parameters,
\[\frac{rc}{2rc -n}=\frac{(t-1)^2(t+2)}{t^3-t^2-4t+6},\]
and it is easy to see that for all $t \geq 3$, we have
\[2 \geq \frac{(t-1)^2(t+2)}{t^3-t^2-4t+6}.\]
Since $r \geq 2$, the third bound in condition (a) is satisfied.
\item[(b)-(f)] These conditions are trivially satisfied given the parametrization of $n$, $r$ and $c$ in terms of $t$ and condition (1)-(4).
\item[(g)] Note that
\begin{align*}
\theta^3-(n-1)&=\frac{1}{8}t^3(t^2-3)^3-\frac{1}{2}(t^2-2)(t^2-1)+1\\
&=\frac{1}{8}(t-1) t^2 (t^2-3) (t (t+1) (t^2-5)+4)
\end{align*}
which is positive for $t\ge 3$.
\item[(h)] If $r=2$, the absolute bound for lines is equivalent to condition (h), as $\overline{m_\theta} = m_\theta$ and $\overline{m_\tau} = m_\tau$.
If $r>2$, then $r\ge 4$. The condition
\[rn\le \frac{1}{2}m_{\theta}(m_{\theta}+1)\]
is equivalent to
\[(r^2-3r+1)m_{\theta}\ge r-1,\]
which holds for all $r\ge 4$. Since $m_{\tau}>m_{\theta}$ for $t\ge 3$, the remaining condition in terms of $m_{\tau}$ is also satisfied.
\item[(i)] We showed that (3) is equivalent to (i) with $\beta = \theta$. If $r \geq 4$ and
\[n>m_{\tau}-r+3,\]
then
\[r<\frac{2(t^4-4t^2+1)}{(t^2-1)(t^2-4)}.\]
But the right hand side is less than four for all $t \geq 3$, hence condition (i) is vacuously satisfied when $\beta = \tau$. 
\end{itemize} 
\ \qed

We show below the first ten feasible parameter sets corresponding to the theorem above. It is not known whether any of these parameter sets correspond to an actual graph.

\begin{table}[here]
\begin{center}
\begin{tabular}{c|c|c|c|c|c|c|c}
$n$ & $r$ & $c$ & $\delta$ & $\theta$ & $\tau$ & $m_\theta$ & $m_\tau$ \\ \hline
 276 & 4 & 56 & 50 & 55 & -5 & 69 & 759 \\
 276 & 16 & 14 & 50 & 55 & -5 & 345 & 3795 \\
 1128 & 6 & 162 & 154 & 161 & -7 & 235 & 5405 \\
 1128 & 54 & 18 & 154 & 161 & -7 & 2491 & 57293 \\
 1128 & 162 & 6 & 154 & 161 & -7 & 7567 & 174041 \\
 1128 & 486 & 2 & 154 & 161 & -7 & 22795 & 524285 \\
 3160 & 4 & 704 & 342 & 351 & -9 & 237 & 9243 \\
 3160 & 8 & 352 & 342 & 351 & -9 & 553 & 21567 \\
 3160 & 64 & 44 & 342 & 351 & -9 & 4977 & 194103 \\
 3160 & 128 & 22 & 342 & 351 & -9 & 10033 & 391287 \\
\end{tabular} \\ \vspace{0.3cm}
Table 1.
\end{center}
\end{table}
Note that if there is an abelian \drackn\ with any of the parameters $(n,r,c)$ above, then there is a regular two-graph with parameters $(n,2,cr/2)$.

\subsection{Case II.b}

\begin{lemma}
Let $X$ be an abelian $(n,r,c)$-\drackn\ such that for some positive integer $t$, we have
\[n=(t^2-1)^2,\quad rc=(t-1)^2(t^2+t-1).\]
Suppose $r\ge 3$. Then $r$ divides $t-1$.
\label{r|t-1}
\end{lemma}
\proof
Suppose that there is a prime $p$ that divides both $r$ and $t^2 + t - 1$. Since $r$ is odd, $p$ is odd. By Corollary \ref{cor:r divides n}, $p$ divides $n$, and $n = (t^2 -1 )^2$. So $p$ divides both $t^2 + t - 1$ and $t^2 -1$, hence $p$ divides $t$, a clear contradiction. Thus, because
\[rc = (t-1)^2 (t^2 + t - 1),\]
it follows that $r $ divides $(t-1)^2$. Hence $r \leq t^2$, and this immediately implies that
\[n > m_\theta -r + 3.\]
By condition (i) of Theorem \ref{feascond}, we now have that $\theta+1$ divides $c$. Since 
\[(t-1)(\theta+1)= r c,\]
this is equivalent to saying that $r$ divides $t-1$.
\qed

The above theorem shows that for $r\ge 3$, the condition that $r$ divides $t-1$ is necessary for Theorem \ref{feascond} to hold. Now we show that this is also sufficient, apart from some lower bounds on $t$, $r$ and $c$.

\begin{theorem}
Suppose that there is a positive integer $t$ defining a parameter set $(n,r,c)$ corresponding to case (II.b) of Theorem \ref{thm:parameters}, and so with $r$ odd. This parameter set satisfies all the feasibility conditions in Theorem \ref{feascond} and Corollary \ref{cor:r divides n} if and only if  $t\ge 3$, $c\ge 2$, $3$ does not divide $r$ and $r$ divides $t-1$.
\end{theorem}
\proof
First we show that the conditions on $t$, $c$ and $r$ are necessary. If $c = 1$, then condition (f) of Theorem \ref{feascond} says that
\[(n-r)^2\le n-1.\]
Note that
\begin{align*}
(n-r)^2-(n-1)&=(n-rc)^2-(n-1)\\
&=(t-1)^4(t+2)^2-(t^2-2)t^2.
\end{align*}
It is easy to see that this is positive for all $t \geq 2$, hence $c > 1$. Moreover, if $t = 2$, then $rc = 5$, and this could only occur with $c = 1$. Hence it also follows that $t \geq 3$.

Now $\theta = (t^2 - 2)t$, so clearly $\theta^3 \neq n-1$. Hence condition (h) of Theorem \ref{feascond} can be applied, and in particular
\[rn \le \frac{1}{2} m_\theta (m_\theta + 1).\]
If $r = 3$, this reduces to $n \leq \sqrt{n}$, a contradiction. If $3$ divides $r$, then we could use Theorem \ref{quotienting} to construct an abelian \drackn\ with the same values of $n$ and $\delta$, hence still meeting the absolute bound with equality, but of index equal to $3$, an absurd. Therefore $3$ cannot divide $r$.

For the converse, we must check that all conditions in Theorem \ref{feascond} are satisfied.

The bounds on condition (a) can be proved similarly to Theorem \ref{thm:caseIb}. Conditions (b) - (d) are trivially satisfied, as well as condition (f). To see (e), note that $n$ is even if and only if $t$ is odd, and since 
\[rc = (t-1)^2(t^2 + t -1)\]
and $r$ divides $t-1$, it follows that $t-1$ divides $c$. So if $n$ is even, then $c$ is even.

As we already mentioned, $\theta^3 > n-1$, hence condition (g) is satisfied. To see (h), note that if $r \geq 4$, then $r < (1/2)(r-1)^2$, hence
\[rn \leq \frac{1}{2} (r-1)^2 n + \frac{1}{2}(r-1) \sqrt{n},\]
which is equivalent to $rn \leq (1/2)m_\theta(m_\theta+1)$. Since $m_\tau > m_\theta$, we have that both inequalities of condition (h) are satisfied.

For condition (i), note that $\theta+1$ divides $c$ if and only if $r$ divides $t-1$, so the case $\beta = \theta$ offers no risk. For the other case, note that
\[\comp{m_{\tau}}=n-\comp{m_{\theta}}=n-\sqrt{n}.\]
Since $r\ge 4$, we have
\begin{align*}
m_{\tau}-r+3-n&=\comp{m_{\tau}}(r-1)-r+3-n\\
&=r(n-\sqrt{n}-1)-(2n-\sqrt{n}-3)\\
&\ge 4(n-\sqrt{n}-1)-(2n-\sqrt{n}-3)\\
&=2n-3\sqrt{n}-1.
\end{align*}
Because $t \geq 3$, we have $n \geq 64$, hence $n < m_\tau -r + 3$. Thus condition (i) is vacuously satisfied in this case.

Finally, Corollary \ref{cor:r divides n} follows from the fact that $r$ divides $t-1$ and thus $r$ divides $n$.
\qed

We show on Table 2 the first ten feasible parameter sets corresponding to the Theorem above. It is not known whether any of these parameter sets corresponds to an actual graph.

\begin{table}[here]
\begin{center}
\begin{tabular}{c|c|c|c|c|c|c|c}
$n$ & $r$ & $c$ & $\delta$ & $\theta$ & $\tau$ & $m_\theta$ & $m_\tau$ \\ \hline
 1225 & 5 & 205 & 198 & 204 & -6 & 140 & 4760 \\
 3969 & 7 & 497 & 488 & 496 & -8 & 378 & 23436 \\
 14400 & 5 & 2620 & 1298 & 1309 & -11 & 480 & 57120 \\
 20449 & 11 & 1705 & 1692 & 1704 & -12 & 1430 & 203060 \\
 38025 & 13 & 2717 & 2702 & 2716 & -14 & 2340 & 453960 \\
 50176 & 7 & 6692 & 3330 & 3345 & -15 & 1344 & 299712 \\
 65025 & 5 & 12195 & 4048 & 4064 & -16 & 1020 & 259080 \\
 104329 & 17 & 5797 & 5778 & 5796 & -18 & 5168 & 1664096 \\
 159201 & 19 & 7961 & 7940 & 7960 & -20 & 7182 & 2858436 \\
 193600 & 5 & 36880 & 9198 & 9219 & -21 & 1760 & 772640 \\
\end{tabular} \\ \vspace{0.3cm}
Table 2.
\end{center}
\end{table}

\section{Final comments and open problems}

Warwick de Launey stated a conjecture in \cite[V.5.18.1]{HandbookDesigns} that there are no generalized Hadamard matrices over non-prime-power order groups. If this conjecture is true, it would imply that the index of all abelian \drackn s with $\delta = -2$ is a prime power. Conversely, if one proves such statement, then De Launey's conjecture will be settled for abelian groups.

A remarkable feature of the feasible parameter sets appearing in Tables 1 and 2 is that $\delta$ increases arbitrarily. Most known \drackn s occur with $\delta \in \{-2,0,2\}$, hence suggesting that these graphs might be very hard to construct.

We are particularly interested in finding a set of equiangular lines satisfying the hypothesis of Theorem \ref{thm:coversfromlines} which yields the construction of a new abelian \drackn. Likewise, it would be very interesting to find a construction for an abelian \drackn\ whose parameters appear in Table 2, thus constructing a set of $d^2$ equiangular lines in $\mathbb{C}^d$.

\section*{Acknowledgements}

This work is supported by C. Godsil's NSERC grant number RGPIN-9439.

\bibliographystyle{amsplain}
\bibliography{el.bib}

\end{document}